\documentclass{article}
\usepackage{amssymb}
\usepackage{amsmath}
\usepackage{enumerate}
\usepackage[latin1]{inputenc}
\def\s{\mathbb{S}}
\def\h{\mathbb{H}}
\def\r{\mathbb{R}}
\def\c{\mathbb{C}}
\def\p{\mathbb{P}}
\newtheorem{remark}{Remark} 
\newtheorem{theorem}{Theorem}
\newtheorem{proposition}{Proposition}
\begin{document}

\title{Minimal Lagrangian submanifolds \\ in the complex hyperbolic space}         
\author{Ildefonso Castro\thanks{Research partially supported by a DGICYT grant No. PB97-0785.}
$\, $  Cristina R. Montealegre \and Francisco Urbano$^*$}        

\date{September 2001}          

\maketitle

\pagestyle{plain}

\begin{abstract}
In this paper we construct new examples of minimal Lagrangian submanifolds in the complex hyperbolic space with large symmetry groups, obtaining three 1-parameter families with cohomogeneity one. We characterize them as the only minimal Lagrangian submanifolds in $\c\h^n$ foliated by umbilical hypersurfaces of Lagrangian subspaces $\r\h^n$ of $\c\h^n$. Several suitable generalizations of the above construction allow us to get new families of minimal  Lagrangian submanifolds in $\c\h^n$ from curves in $\c\h^1$ and $(n-1)$-dimensional minimal  Lagrangian submanifolds of the complex space forms $\c\p^{n-1}$, $\c\h^{n-1}$ and $\c^{n-1}$. Similar constructions are made in the complex projective space $\c\p^n$.
\end{abstract}

\section{Introduction}       

Special Lagrangian submanifolds of complex Euclidean space $\c^n$ (or of a Calabi-Yau manifold) have been studied widely over the last few years. They have appeared in Mathematical Physics in [SYZ], where A. Strominger, S.T. Yau and E. Zaslow proposed an explanation of mirror symmetry of a Calabi-Yau manifold in terms of the moduli space of special Lagrangian submanifolds. These submanifolds are volume minimizing and, in particular, they are minimal submanifolds. Furthermore any oriented minimal Lagrangian submanifold of $\c^n$ (or of a Calabi-Yau manifold) is a special Lagrangian submanifold with respect to one of the 1-parameter family of special Lagrangian calibrations which this kind of Kaehler manifolds has. (see [HL, Proposition 2.17]).

A very important problem here is finding (non-trivial) examples of special Lagrangian submanifiolds (i.e. oriented minimal Lagrangian submanifolds). In [HL] R. Harvey and H.B. Lawson constructed the first examples in $\c^n$, where we point out the Lagrangian catenoid one ([HL], example III.3.B), and more recently D.D. Joyce ([J1,J2,J3,J4]) and M. Haskins ([H]) have developed methods for constructing important families of special Lagrangian submanifolds of $\c^n$. We are particularly interested in the examples with large symmetry groups ([J1]), i.e. invariant under the action of certain subgroups of the isometries group of $\c^n$.

Following some ideas in the mentioned papers, we construct examples of minimal Lagrangian submanifolds of the complex hyperbolic space $\c\h^n$ with large symmetry groups. In particular, we consider the groups of isometries of the sphere $\s^{n-1}$, the real hyperbolic space $\r\h^{n-1}$ and Euclidean space $\r^{n-1}$:  $SO(n)$, $SO^1_0(n)$ and $SO(n\!-\!1)\propto \r^{n-1}$ respectively, acting on $\c\h^n$ as holomorphic isometries (see section 2.1). In Theorems 1, 2 and 3 we classify respectively the minimal Lagrangian submanifolds of $\c\h^n$ invariant under the groups $SO(n)$, $SO^1_0(n)$ and $SO(n\!-\!1)\propto\r^{n-1}$ acting on $\c\h^n$ as holomorphic isometries. In each result we obtain a 1-parameter family of minimal Lagrangian submanifolds $M$ in $\c\h^n$ with cohomogeneity one, i.e. the orbits of the symmetry group are of codimension one in $M$. In particular $M$ is foliated by a  1-parameter family of orbits parameterized by $s\in \r$, which are geodesic spheres, tubes over hyperplanes and horospheres respectively (i.e. umbilical hypersurfaces) of Lagrangian subspaces $\r\h^n_s$ of $\c\h^n$. In Theorem 4 we characterize the above examples as the only minimal Lagrangian submanifolds of $\c\h^n$ foliated by umbilical hypersurfaces of Lagrangian subspaces $\r\h^n$ of $\c\h^n$.
In [CU2] a similar result to Theorem 4 characterizing the Lagrangian catenoid in $\c^n$ was proved.

Following an idea given independently in [H, Theorem A], [J1, Theorem 6.4] and [CU2, Remark 1] (see Remark 1 in this paper for a better understanding), in Propositions 3 and 4 we construct families of minimal Lagrangian submanifolds of $\c\h^n$ from curves in $\c\h^1$ and $(n\!-\!1)$-dimensional minimal Lagrangian submanifolds of the complex space forms $\c\p^{n-1}$, $\c\h^{n-1}$ and $\c^{n-1}$. The examples described in Theorems 1, 2 and 3 are the simplest ones in the above construction.

Similar results to the ones given in $\c\h^n$ can be obtained in the complex projective space $\c\p^n$ (see Theorem 5 and Proposition 6). In this case, less families of minimal Lagrangian submanifolds appear, because there is only one family of umbilical hypersurfaces of the Lagrangian subspaces $\r\p^n$ of $\c\p^n$: the geodesic spheres. So, we study the case of
$\c\h^n$ in this paper as the descriptions and the proofs are far more difficult.

\section{Preliminaries}

\subsection{The complex hyperbolic space}

In this paper we will consider the following model for the complex hyperbolic space.
In $\c^{n+1}$ we take the Hermitian form $(,)$ given by
\[
(z,w)=\sum_{i=1}^n z_i\bar{w}_i-z_{n+1}\bar{w}_{n+1},
\]
for $z,w\in\c^{n+1}$, where $\bar z$ stands for the conjugate of $z$. If
\[
\h^{2n+1}_1=\{z\in\c^{n+1}\,/\,(z,z)=-1\}
\]
is the anti-De Sitter space, then $\Re\,(,)$ induces on $\h^{2n+1}_1$ a Lorentzian metric of constant curvature $-1$ (Here $\Re$ means real part). If $(\c\h^n=\h^{2n+1}_1/\s^1,\langle,\rangle)$ denotes the complex hyperbolic space of constant holomorphic sectional curvature $-4$, then
\[
\c\h^n=\{\Pi(z)=[z]\,/\,z=(z_1,\dots,z_{n+1})\in\h^{2n+1}_1\},
\]
where $\Pi:\h^{2n+1}_1\rightarrow\c\h^n$ is the Hopf projection. The metric $\Re\,(,)$ becomes  $\Pi$ in a pseudo-Riemannian submersion. The complex structure of $\c^{n+1}$ induces, via $\Pi$, the canonical complex structure $J$ on $\c\h^n$. The K\"{a}hler two--form $\Omega$ in $\c\h^n$ is defined by $\Omega(u,v)=\langle Ju,v\rangle$. 
We recall that $\c\h^n$ has a smooth compactification $\c\h^n\cup\s^{2n-1}(\infty)$, where $\s^{2n-1}(\infty)=\pi({\cal N})$,
\[
{\cal N}=\{z\in\c^{n+1}-\{0\}\,/\,(z,z)=0\},
\]
and $\pi:{\cal N}\rightarrow\s^{2n-1}(\infty)$ is the projection given by the natural action of $\c^*$ over ${\cal N}$.

Moreover, in the paper we will denote by $\c\p^n$ the complex $n$-dimensional complex projective space endowed with the Fubini-Study metric of constant holomorphic sectional curvature 4, and by $\Pi:\s^{2n+1}\rightarrow\c\p^n$ the Hopf fibration from the $(2n\!+\!1)$-dimensional unit sphere $\s^{2n+1}$. We also denote  the complex structure and the K\"{a}hler two-form in $\c\p^n$ by $J$ and $\Omega$ respectively.

If $U^1(n\!+\!1)$ is the group preserving the Hermitian form $(,)$, then
\[
U^1(n\!+\!1)=\{A\in GL(n\!+\!1,\c)\,/\,\bar A^tSA=S\}
\] 
where
\[
S=\left( \begin{array}{c|c}

\mbox{\large I}_n &  \\
\hline
&  -1
\end{array}\right) ,
\]
with $I_n$ the identity matrix of order $n$.

Then $PU^1(n+1)=U^1(n+1)/\s^1$ is the group of the holomorphic isometries of $(\c\h^n,\langle,\rangle)$. 

\vspace{0.3cm}

Along the paper we will work with the special orthogonal group $SO(n)$, the identity component of the indefinite special orthogonal group $SO^1_0(n)$ and the group of isometries of Euclidean $(n \! - \! 1)$-space $SO(n \! - \! 1)\propto \r^{n-1}$, acting on $\c\h^n$ as subgroups of holomorphic isometries, in the following ways:
\[
A\in SO(n)\longmapsto \left[\left(\begin{array}{c|c}

A &  \\
\hline
&  1
\end{array}\right)\right]\in PU^1(n+1),
\]
\[
A\in SO^1_0(n)\longmapsto \left[\left(\begin{array}{c|c}
1 & \\ 
\hline
&  A   

\end{array}\right)\right]\in PU^1(n+1),
\]
\[
(A,a)\in SO(n-1)\propto \r^{n-1}\longmapsto \left[\left(
\begin{array}{c|c|c}
A & Aa^t  & Aa^t \\
\hline
-a & 1-\frac{|a|^2}{2} & -\frac{|a|^2}{2} \\
\hline
a & \frac{|a|^2}{2} & 1+\frac{|a|^2}{2}
\end{array} \right)\right]\in PU^1(n+1),
\]
where $a=\left( a_1,\dots,a_{n-1}\right)$.
Here $[\;\;]$ stands for class in $U^1(n+1)/\s^1$.

\subsection{Lagrangian submanifolds in $\c\h^n$}

Let $\phi$ be an isometric immersion of a Riemannian $n$-manifold
$M$ in $\c\h^n$ (resp. $\c\p^n$). $\phi $ is called {\em Lagrangian} if $\phi^*\Omega\equiv 0$.
We denote the Levi-Civita connection of $M$ and the connection on the
normal bundle by $\nabla $ and $\nabla ^\perp $ respectively.
The second fundamental form will be denoted by $\sigma $.
If $\phi $ is Lagrangian, the formulas of Gauss and Weingarten  lead to
\begin{eqnarray*}
\nabla^\perp_X JY = J \nabla _X Y , 
\end{eqnarray*}
and the trilinear form $\langle \sigma (X,Y) , JZ \rangle$
is totally symmetric for any tangent vector fields $X$, $Y$ and $Z$.

If $\phi:M\longrightarrow \c\h^{n}$ (resp. $\c\p^n$) is a Lagrangian immersion of a simply-connected manifold $M$, then $\phi$ has a horizontal lift with respect to the Hopf fibration to $\h^{2n+1}_1$ (resp. $\s^{2n+1}$), which is unique up to isometries. We will denote  this horizontal lift by $\tilde{\phi}$.
We note that only Lagrangian immersions in $\c\h^n$ (resp. $\c\p^n$) have (locally) horizontal lifts. Horizontal immersions from n-manifolds in $\h^{2n+1}_1$ (resp. $\s^{2n+1}$) are called {\em Legendrian immersions} (see [H]). So we can paraphrase the above reasoning as follows: {\em Lagrangian immersions in $\c\h^n$ (resp. $\c\p^n$) are locally projections of Legendrian immersions in $\h^{2n+1}_1$} (resp. $\s^{2n+1}$).

If $H$ is the mean curvature vector of the immersion $\phi:M\longrightarrow \c\h^{n}$, then $\phi$ is  called minimal if $H=0$. The minimality means that the submanifold is critical for compact supported variations of the volume functional. In [O], the second variation of the volume functional was studied for minimal Lagrangian submanifolds of Kaehler manifolds. Among other things, it was proved that {\em minimal Lagrangian submanifolds in $\c\h^n$ are stable and without nullity.}

\vspace{0.3cm}

Let $\langle\langle,\rangle\rangle$ be the restriction of $(,)$ to $\r^{n+1}\equiv\Re \c^{n+1}$. The real hyperbolic space $\r\h^n$ endowed with its canonical metric of constant sectional curvature $-1$ is defined as the following hypersurface of $(\r^{n+1},\langle\langle,\rangle\rangle)$:
$$\r\h^n=\{x\in\r^{n+1}\,/\,\langle\langle x,x\rangle\rangle=-1,\,x_{n+1}\geq 1\}.$$
We recall that $\r\h^n$ has also a smooth compactification $\r\h^n\cup\s^{n-1}(\infty)$, where $\s^{n-1}(\infty)=\pi({\cal N})$, with ${\cal N}=\{x\in\r^{n+1}-\{0\}\,/\,\langle\langle x,x\rangle\rangle=0\}$ the light cone and $\pi$ the projection given by the natural action of $\r^*$ over ${\cal N}$. In addition, $SO^1_0(n+1)$ is a group of isometries of $(\r\h^n,\langle \langle,\rangle\rangle )$.

$\r\h^n$ can be isometrically embedded in $\c\h^n$ as a totally geodesic Lagrangian submanifold in the standard way
\[
x\in\r\h^n \mapsto [x]\in\c\h^n.
\]
Moreover, up to congruences, it is the only totally geodesic Lagrangian submanifold of $\c\h^n$.
It is also interesting to point out (for later  use in section 4) that the totally umbilical submanifolds of $\c\h^n$ (which were classified in [ChO]) are either totally geodesic or umbilical submanifolds of totally geodesic Lagrangian submanifolds. So, {\em up to congruences, the (n-1)-dimensional totally umbilical (non totally geodesic) submanifolds of $\c\h^n$ are the umbilical hypersurfaces of $\r\h^n$ embedded in $\c\h^n$ in the above way}. Up to congruences, the umbilical hypersurfaces of $\r\h^n$ can be described in the following way:
\begin{enumerate}
\item  {\em Geodesic spheres}. Given $r>0$, let $\psi:\s^{n-1}\rightarrow \r\h^n$ be the embedding given by
\[ 
\psi(x)=\left (\sinh r \,\, x,\cosh r\right ).
\]
Then $\psi(\s^{n-1})$ is the geodesic sphere of $\r\h^n$ of center $(0,\dots,0,1)$ and radius $r$.
\item  {\em Tubes over hyperplanes}. Given $r>0$, let $\psi:\r\h^{n-1}\rightarrow \r\h^n$ be the embedding given by
\[ 
\psi(x)=\left (\sinh r,\cosh r \,\, x\right ).
\]
Then $\psi(\r\h^{n-1})$ is the tube to distance $r$ over the hyperplane dual to $(1,0,\dots,0)$.
\item {\em Horosphere}. Let $\psi:\r^{n-1}\rightarrow \r\h^n$ be the embedding given by
\[ 
\psi(x)=\left (x,\frac{|x|^2}{2},\frac{|x|^2}{2}+1\right ).
\]
Then $\psi(\r^{n-1})$ is a horosphere of $\r\h^n$ with infinity point $\pi(0,\dots,0,1,1)$.
\end{enumerate} 
We will refer to these examples as {\em (n-1)-geodesic spheres}, {\em (n-1)-tubes over hyperplanes} and {\em (n-1)-horospheres} of $\c\h^n$.
 
\section{Examples of minimal Lagrangian submanifolds with symmetries}

In this section we are going to describe the minimal Lagrangian submanifolds of $\c\h^n$ invariant under the actions of $SO(n)$, $SO_0^1(n)$ and $SO({n\!-\!1})\propto\r^{n-1}$ as subgroups of isometries of $\c\h^n$ given in \S2.1. These examples of minimal Lagrangian submanifolds can be regarded as the simplest ones.

\subsection{Examples invariant under $SO(n)$}
\begin{theorem}
For any $\rho > 0$, there exists a minimal (non totally geodesic) Lagrangian embedding
\[ \Phi_{\rho} : \r \times \s^{n-1} \longrightarrow \c\h^n \]
defined by
\[ 
\Phi_{\rho} (s,x) = \left[\left( 
\sinh r(s) \, e^{i \int_0^s \frac{dt}{\sinh^{n+1} r(t)} } \, x  \, , \, 
\cosh r(s) \, e^{i \int_0^s \frac{\tanh^2 r(t) \, dt} {\sinh^{n+1} r(t)} }
\right)\right],
\]
where $r(s)$, $s\in \r$, is the only solution to
\begin{equation}
r''\sinh r \cosh r  = (1-(r')^2)(\sinh^2 r + n\cosh^2 r), \, r(0)=\rho, \,
r'(0)=0.
\end{equation}
 
$\Phi_\rho$ is
invariant under the action of $SO(n)$ and satisfies
\[
\int_{\r\times\s^{n-1}}|\sigma|^n dv\,<\infty,
\]
where $dv$ is the canonical measure of the complete induced metric $ds^2+\sinh^2r(s)g_0$, with $g_0$ the canonical metric of the unit sphere $\s^{n-1}$. 

Moreover, any minimal (non totally geodesic) Lagrangian immersion in $\c\h^n$ invariant under the action of $SO(n)$ is congruent to an open subset of some of the above submanifolds.
\end{theorem}

{\it Proof:\/}
We start the proof  with the analysis of the differential equation (1). 
The energy integral of (1) is given by
\[
(1-(r')^2)\cosh^2 r \sinh^{2n} r = {\rm constant},
\]
which is equivalent to
\[
(r')^2 + \frac{ch^2_{\rho}sh^{2n}_{\rho}}{\cosh^2 r \sinh^{2n} r}=1,
\]
where $ch_{\rho} = \cosh \rho$ and $sh_{\rho} = \sinh \rho$.

Using qualitative theory of o.d.e., the existence  and uniqueness of solution $r(s)$ to (1)
is guaranteed, defined on the whole $\r$, for any initial condition $\rho=r(0)>0$. This will be the only absolute minimum of $r(s)$ and $r(-s)=r(s), \forall s\in \r$. 

It is easy to prove that $\Phi_\rho $ is a Lagrangian immersion invariant under the action of $SO(n)$.
The induced metric $ds^2+\sinh^2r(s)g_0$ is a complete metric, since  $\sinh^2r(s) \geq \sinh^2\rho$. We consider the orthonormal frame for this metric given by
\[
e_1=\partial_s, \, e_j=v_j / \sinh r, \, j=2,\dots,n ,
\]
where $\{ v_2,\dots,v_n \}$ is an orthonormal frame of $(\s^{n-1},g_0)$.
Then we can compute the second fundamental form $\sigma $ of $\Phi_{\rho}$:
\[
\sigma(e_1,e_1)=-\frac{(n-1)Je_1}{\sinh^{n+1}r}, \quad
\sigma(e_1,e_j)=\frac{Je_j}{\sinh^{n+1}r}, \quad
\sigma(e_j,e_k)=\frac{\delta_{jk}Je_1}{\sinh^{n+1}r} .
\]
Using the information above, it is not difficult to prove that 
$\Phi_{\rho}$ is minimal and, after some computations, we also arrive at
\[
\int_{\r\times\s^{n-1}}|\sigma|^n dv = 2((n+2)(n-1))^{n/2}c_{n-1} 
\int_{0}^{+\infty} \frac{ds}{\sinh^{n^2+1}r(s)},
\]
where $c_{n-1}$ denotes the volume of $(\s^{n-1},g_0)$. Making the change of variable $t=\sinh r(s)$, we get:
\[ \int_{0}^{+\infty} \frac{ds}{\sinh^{n^2+1}r(s)} =
\int_{\sinh \rho}^{+\infty} \frac{dt}{t^{n^2-n+1}\sqrt{t^{2n+2}+t^{2n}-ch^2_{\rho}sh^{2n}_{\rho}}}.
\]
The last is an hyperelliptyc integral and we can prove that converges using
numerical methods.

Now we prove that $\Phi_\rho$ is an embedding. Suppose $\Phi_\rho(s,x)=
\Phi_\rho(\hat s,\hat x)$.
This implies there exists $\theta \in \r$ such that the horizontal lift $\tilde{\Phi}_\rho$ of our immersion verifies:
\[
\tilde{\Phi}_\rho(\hat s,\hat{x})=e^{i\theta}\tilde{\Phi}_\rho(s,x).
\] 
From the definition of $\Phi_\rho$ we deduce that $r(\hat s)=r(s)$ and so
$\hat s=\pm s$. If $\hat s = s$, necessarily $\hat x = x$. But if $\hat s=-s$, we get 
\[
\hat{x}= e^{2i\int_0^s \frac{dt}{\cosh^2 r(t) \sinh^{n+1}r(t)}}x .
\]
Using a similar reasoning as above, we can check that the increasing function $s \rightarrow 2\int_0^s \frac{dt}{\cosh^2 r(t) \sinh^{n+1}r(t)}$ never reaches the value $\pi$. Since the coordinates of $x$ and $\hat{x}$ are real numbers it is impossible that $\hat s=-s$ and hence $\Phi_\rho$ must be an embedding.

\vspace{0.2cm}

Conversely, let $\phi:M\longrightarrow \c\h^n$ be a non totally geodesic minimal Lagrangian immersion invariant under the action of $SO(n)$, and $\tilde{\phi}$ a local horizontal lift of $\phi$ to $\h^{2n+1}_1$. Let $p$ any point of $M$ and $z=(z_1,\dots,z_{n+1})=\tilde{\phi}(p)$. As $\phi$ is invariant under the action of $SO(n)$, for any matrix $A$ in the Lie algebra of $SO(n)$, the curve $s\rightarrow [ze^{s\hat A}]$ with 
\[
\hat A=\left( \begin{array}{c|c}

\mbox{ A} &  \\
\hline
&  0
\end{array}\right),
\]
lies in the submanifold, and then its tangent vector at $s=0$ satisfies
$$\Pi_*(z\hat A+(z\hat A,z)z)\in \phi_*(T_pM).$$
Using that $\phi$ is Lagrangian, this implies that
\[
\Im (z\hat A\hat B\bar z^t)=0,
\]
for any $n$-matrixes $A,B$ in the Lie algebra of $SO(n)$. If $n\geq 3$, from here it is easy to see that 
$\Re (z_1,\dots,z_n)$ and $\Im (z_1,\dots,z_n)$ are linearly dependent. As $SO(n)$ acts transitively on $\s^{n-1}$, we obtain that $z$ is in the orbit (under the action of $SO(n)$ described above) of the point $(a+ib,0,\dots,0,z_{n+1})$, with $a^2+b^2=|z_{n+1}|^2-1$. This reasoning implies that locally $\tilde{\phi}$ is the orbit under the action of $SO(n)$ of a curve in $\h^3_1\equiv\h^{2n+1}_1\cap\{z_2=\dots=z_n=0\}$. So $M$ is locally $I\times\s^{n-1}$, with $I$ an interval in $\r$; moreover, the lift $\tilde{\phi}:I\times\s^{n-1}\rightarrow\h^{2n+1}_1$ is given by
\[
\tilde{\phi}(s,x)=(\gamma_1(s)x,\gamma_2(s)),
\]
where $\gamma(s)=(\gamma_1(s),\gamma_2(s))$ is an horizontal curve in $\h^3_1$. If $n=2$ it is also easy to get the above expression. The horizontality of the curve $\gamma$ in $\h^3_1$ means that we can find real functions $r=r(s)>0$ and $f=f(s)$, such that
\[
\gamma(s)=\left (\sinh r(s)e^{i\int_{s_0}^s f(t)dt} \, , \,
\cosh r(s)e^{i\int_{s_0}^s f(t)\tanh^2r(t) dt}\right ),
\]
with $s_0\in I$. Now, using that $\phi$ is a minimal immersion we are going to determine the functions $f$ and $r$. After a long but straightforward computation, one can prove that the immersion $\phi$ is minimal if and only if $f$ and $r$ satisfy the following equation
\begin{equation}
f r''\tanh r=f^3\tanh^2r\,(n+\tanh^2r)+(n+1)(r')^2 f+f' r'\tanh r.
\end{equation}
If $r$ is constant, necessarily $f\equiv0$ and $\gamma$ degenerates into a point. 

In order to analyze the equation (2) in the non trivial case, we assume that $\gamma$ is parameterized by the arc, i.e., $|\gamma'|=1$. By computing $|\gamma'|$, we get:
\[
(r')^2+f^2\tanh^2r=1.
\]
Deriving this equation and using it again, the equation (2) becomes in
\[
(n+1)fr'+f'\tanh r=0.
\]
The solution $f\equiv 0$ says that $r(s)$ is a linear map and this leads to the totally geodesic case.
The general solution to the above equation is
\[
f(s)=\frac{a}{\sinh^{n+1}r(s)},\quad a>0.
\]
In this way, we have proved that $r(s)$ must satisfy the equation
\[
(r')^2 + \frac{a^2}{\cosh^2 r \sinh^{2n} r}=1.
\]
The solutions to this differential equation are defined in the whole $\r$ and have only one critical point. Therefore we can take $s_0=0$ in the definition of $\gamma$ and consider $r'(0)=0$ up to a translation of parameter. So $a^2=\cosh^2 r(0) \sinh^{2n} r(0)$ and then $r$ is a solution to (1).$_{\diamondsuit}$

\vspace{1cm}

We observe that for each $s\in\r$, $\Phi_{\rho}(\{s\}\times\s^{n-1})$ is a geodesic sphere of radius $r(s)$ and center $[(0,\dots,0,1)]$ in the Lagrangian subspace $\r\h^n_s$ of $\c\h^n$ defined by
\[
\r\h^n_s=\{[(x_1,\dots,x_{n+1})A(s)]\,/\,x_i\in\r\,,\,\sum_{i=1}^nx_i^2-x_{n+1}^2=-1\,,\,x_{n+1}\geq 1\}
\]
where $A(s)$ is the matrix of $U^1(n+1)$ defined by
\[
\ A(s)=\left( \begin{array}{c|c}

e^{ia(s)}I_{n} &  \\
\hline
&  e^{ib(s)}
\end{array}\right),
\]
being $a(s)=\int_0^s \frac{dt}{\sinh^{n+1} r(t)}$ and $b(s)=\int_0^s \frac{\tanh^2 r(t) \, dt} {\sinh^{n+1} r(t)}$. Moreover, if $s\not= s'$ then $\r\h^n_s\cap\r\h^n_{s'}=[(0,\dots,0,1)]$. Hence $\{\Phi(\{s\}\times\s^{n-1})\,,s\in\r\,\}$ defines a foliation on the minimal Lagrangian submanifold by $(n-1)$-geodesic spheres of $\c\h^n$.

In a more general context we can classify pairs of Lagrangian subspaces of $\c\h^n$ intersecting only in a point. (Compare with Proposition 6.2 in [J1]). 

\begin{proposition}
Let $\r\h^n_a$ and $\r\h^n_b$ two Lagrangian subspaces of $\c\h^n$ which intersect only at $[(0,\dots,0,1)]$. Then there exist $\theta_1,\dots,\theta_n\in (0,\pi)$ and $A\in U^1(n\!+\!1)$ such that
\[
\r\h^n_a=\{[(x_1,\dots,x_{n+1})A]\,/\,x_i\in\r\,,\,\sum_{i=1}^nx_i^2-x_{n+1}^2=-1\,,\,x_{n+1}\geq 1\}
\]
and 
\[
\r\h^n_b=\{[(e^{i\theta_1}x_1,\dots,e^{i\theta_n}x_n,x_{n+1})A]\,/\,x_i\in\r\,,\,\sum_{i=1}^nx_i^2-x_{n+1}^2=-1\,,\,x_{n+1}\geq 1\}.
\]
\end{proposition}

Two Lagrangian subspaces $\r\h^n_a$ and $\r\h^n_b$ which intersect only at $[(0,\dots,0,1)]$ with $\theta_1=\dots =\theta_n$ are said to be in {\em ``normal'' position}. In particular, in our family of Lagrangian subspaces $\{ \r\h^n_s\,/\,s\in\r\}$ every two Lagrangian subspaces $\r\h^n_s$ and $\r\h^n_{s'}$ are in normal position.
\begin{proposition}
Let $\phi:M\rightarrow\c\h^n$ be a minimal Lagrangian immersion of a compact manifold with boundary $\partial M$. If $\phi(\partial M)$ is the union of two geodesic spheres centered at $[(0,\dots,0,1)]$ in two Lagrangian subspaces in normal position, then $\phi$ is congruent to some of the examples given in Theorem 1.
\end{proposition}
{\it Proof: \/}
It is clear that, up to a holomorphic isometry of $\c\h^n$, the Lagrangian subspaces in normal position can be taken as
\[
\r\h^n_1=\{[(x_1,\dots,x_{n+1})]\,/\,x_i\in\r\,,\,\sum_{i=1}^nx_i^2-x_{n+1}^2=-1\,,\,x_{n+1}\geq 1\}
\]
and 
\[
\r\h^n_2=\{[(e^{i\theta}x_1,\dots,e^{i\theta}x_n,x_{n+1})]\,/\,x_i\in\r\,,\,\sum_{i=1}^nx_i^2-x_{n+1}^2=-1\,,\,x_{n+1}\geq 1\}.
\]
Now, these Lagrangian subspaces and their corresponding geodesic spheres centered at $[(0,\dots,0,1)]$ are invariant under the action of the group $SO(n)$ on $\c\h^n$ (see section 2). Hence if $X$ is a Killing vector field in the Lie algebra of $SO(n)$, then its restriction to the submanifold is a Jacobi field on $M$ vanishing on $\partial M$. As the nullity of the submanifold is zero, X also vanishes along the submanifold $M$. This means that the submanifold is invariant under the action of $SO(n)$ and then the result follows from Theorem 1.$_{\diamondsuit}$

\subsection{Examples invariant under $SO^1_0(n)$}
\begin{theorem}
For any $\rho> 0$, there exists a minimal (non totally geodesic) Lagrangian embedding
\[ \Psi_{\rho} : \r \times \r\h^{n-1} \longrightarrow \c\h^n \]
defined by
\[ 
\Psi_{\rho} (s,x) = \left[\left( 
\sinh r(s) \, e^{i \int_0^s \frac{\coth^2 r(t) \, dt}{\cosh^{n+1} r(t)} } \,  ,  \,
\cosh r(s) \, e^{i \int_0^s \frac{dt} {\cosh^{n+1} r(t)} } \, x
\right)\right],
\]
where $r(s)$, $s\in \r$, is the only solution to
\begin{equation}
r''\sinh r \cosh r  = (1-(r')^2)(\cosh^2 r + n\sinh^2 r), \, r(0)=\rho, \,
r'(0)=0.
\end{equation}

$\Psi_\rho$ is
invariant under the action of $SO^1_0(n)$ and satisfies
\[
\int_{\r\times\r\h^{n-1}}|\sigma|^n dv\,<\infty,
\]
where $dv$ is the canonical measure of the complete induced metric $ds^2+\cosh^2r(s)$ $ \langle \langle, \rangle \rangle$.

Moreover, any minimal (non totally geodesic) Lagrangian immersion in $\c\h^n$ invariant under the action of $SO^1_0(n)$ is congruent to an open subset of some of the above submanifolds.
\end{theorem}

As the proof of Theorem 2 is similar to the proof of Theorem 1, it will be omitted.

\subsection{Examples invariant under $SO(n\hspace{-3pt}-\hspace{-3pt}1)\propto\r^{n-1}$}

\begin{theorem}
For any $\rho > 0$, there exists a minimal Lagrangian embedding
\[ \Upsilon_{\rho} : \r \times \r^{n-1} \longrightarrow \c\h^n \]
defined by
\[
\begin{array}{l}
\Upsilon_{\rho} (s,x)    = \\
 \left[e^{i A_{n+1}(s)} \left(
r(s)x,
\frac{1+r(s)^2(|x|^2-1-2i A_{n+3}(s))}{2r(s)},
\frac{1+r(s)^2(|x|^2+1-2i A_{n+3}(s))}{2r(s)}
\right) \right],
\end{array}
\]
where $A_n(s)=\int_0^s dt/r(t)^n $ and
\[
r(s)=\rho \cosh^{\frac{1}{n+1}} ((n+1)s).
\]

$\Upsilon_\rho$ is invariant under the action of $SO(n\hspace{-3pt}-\hspace{-3pt}1)\propto\r^{n-1}$ and satisfies
\[
\int_{\r\times\r^{n-1}}|\sigma|^ndv\,<\infty,
\]
where $dv$ is the canonical measure of the complete induced metric $ds^2+r(s)^2 \langle, \rangle$, where $\langle, \rangle $ is the canonical metric of Euclidean space $\r^{n-1}$. 

Moreover, any minimal (non totally geodesic) Lagrangian immersion in $\c\h^n$ invariant under the action of $SO(n\hspace{-3pt}-\hspace{-3pt}1)\propto\r^{n-1}$ is congruent to an open subset of some of the above submanifolds.
\end{theorem}

{\it Proof: \/} In a similar way that in \S3.1, the geometric properties of $\Upsilon_{\rho} $ can be checked now from the explicit expressions given in the statement of the Theorem.

\vspace{0.2cm}

Conversely, let $\phi:M\rightarrow \c\h^n$ be a non totally geodesic minimal Lagrangian immersion invariant under the action of $SO(n\hspace{-3pt}-\hspace{-3pt}1)\propto\r^{n-1}$, and $\tilde{\phi}$ a local horizontal lift of $\phi$ to $\h^{2n+1}_1$. Let $p$ a point of $M$ and $z=(z_1,\dots,z_{n+1})=\tilde{\phi}(p)$. As $\phi$ is invariant under the action of $SO(n\hspace{-3pt}-\hspace{-3pt}1)\propto\r^{n-1}$ then for any $(A,a)$ in the Lie algebra of $SO(n\hspace{-3pt}-\hspace{-3pt}1)\propto\r^{n-1}$, the curve $s\rightarrow [ze^{s\hat A}]$ with 
\[
\hat A=\left(
\begin{array}{c|c|c}
A & a^t  & a^t \\
\hline
-a & 0 & 0 \\
\hline
a & 0 & 0
\end{array} \right)
\]
lies in the submanifold, and then 
$$\pi_*(z\hat A+(z\hat A,z)z)\in \phi_*(T_pM).$$
Since $\phi$ is Lagrangian, we can deduce that
\[
\Im (z\hat A\hat B^t\bar z^t +(z\hat A,z)(z,z\hat B))=0,
\]
for any $\hat A,\hat B$ in the Lie algebra of $SO(n\hspace{-3pt}-\hspace{-3pt}1)\propto\r^{n-1}$. If $n\geq 3$, from here it is easy to see that 
$(z_1,\dots,z_{n-1})=(z_{n+1}\hspace{-3pt}-\hspace{-3pt}z_n)(x_1,\dots,x_{n-1})$, with $x=(x_1,\dots,x_{n-1})\in\r^{n-1}$. As $SO(n\hspace{-3pt}-\hspace{-3pt}1)\propto\r^{n-1}$ acts transitively on $\r^{n-1}$, we obtain that $z$ is in the orbit under the action of $SO(n\hspace{-3pt}-\hspace{-3pt}1)\propto\r^{n-1}$ described above, of the point 
\[
\left (0,\dots,0,z_n-(z_{n+1}\hspace{-3pt}-\hspace{-3pt}z_n)\frac{|x|^2}{2},z_{n+1}-(z_{n+1}\hspace{-3pt}-\hspace{-3pt}z_n)\frac{|x|^2}{2}\right ).
\]
This reasoning implies that locally $\tilde{\phi}$ is the orbit under the action of $SO(n\hspace{-3pt}-\hspace{-3pt}1)\propto\r^{n-1}$ of a curve in $\h^3_1\equiv\h^{2n+1}_1\cap\{z_1=\dots=z_{n-1}=0\}$. So $M$ is locally $I\times\r^{n-1}$, with $I$ an interval in $\r$ and the lift $\tilde{\phi}:I\times\r^{n-1}\rightarrow\h^{2n+1}_1$ is given by
\[
\tilde{\phi}(s,x)=\left (\gamma_2(s)\hspace{-3pt}-\hspace{-3pt}\gamma_1(s)\right )\left (x,\frac{|x|^2}{2},\frac{|x|^2}{2}\right )+\left (0,\gamma_1(s),\gamma_2(s)\right )\]
where $\gamma(s)=(\gamma_1(s),\gamma_2(s))$ is an horizontal curve in $\h^3_1$. If $n=2$ it is also easy to get the above expression. 
Writing $(\gamma_2-\gamma_1)(s)=r(s)e^{i\int_{s_0}^sf(t)dt}$ for certain real functions $r=r(s)>0$ and $f=f(s)$, the horizontality of $\gamma$ implies that
\[
\gamma(s)=e^{i\int_{s_0}^sf(t)dt}\left (\frac{1-r(s)^2}{2r(s)}-ir(s)\!\int_{s_0}^s\frac{f(t)}{r(t)^2}dt \, , \, \frac{1+r(s)^2}{2r(s)}-ir(s)\!\int_{s_0}^s\frac{f(t)}{r(t)^2}dt\right ),
\]
with $s_0\in I$.
A similar reasoning to that in \S3.1 translates the minimallity of the immersion $\phi$  into 
\begin{equation}
(n+1)f((r')^2+r^2f^2)-f r r''+f' r r'+f(r')^2=0.
\end{equation}
If $r$ is constant, necessarily $f\equiv0$ and $\gamma$ degenerates into a point. 

To analyze the equation (4), we assume that $\gamma$ is parameterized by the arc, i.e., $|\gamma'|=1$. By computing $|\gamma'|$ we get
\[
(r'/r)^2+f^2=1.
\]
Deriving this equation and using it again in (4), one get
\[
(n+1)fr'+f'r=0.
\]
The solution $f\equiv 0$ says that $r(s)=\mu e^{\pm s}$ and the immersion $\phi$ is totally geodesic. Otherwise the general solution is given by
\[
f(s)=\frac{a}{r(s)^{n+1}},\quad a>0,
\]
and so we have proved that $r(s)$ must satisfy the equation
\[
(r')^2 + \frac{a^2}{r^{2n}}=r^2,
\]
whose general solution, up to a translation of parameter, is given in the statement of the Theorem putting $a=\rho^{n+1}$.
$_{\diamondsuit}$

\section{More examples of minimal Lagrangian submanifolds}

The examples given in Theorems 1, 2 and 3 have a particular common way to be constructed and the analysis of this construction will be the key to give new examples of minimal Lagrangian submanifolds in $\c\h^n$. In fact, the examples in Theorems 1, 2 and 3 are constructed, respectively, in the following way:
\[
(s,x)\in\r\times\s^{n-1}\mapsto[(\gamma_1(s)\,x,\gamma_2(s))]\in\c\h^n,
\]
\[
(s,x)\in\r\times\c\h^{n-1}\mapsto[(\gamma_1(s),\gamma_2(s)\,x)]\in\c\h^n,
\]
\[
(s,x)\in\r\times\r^{n-1}\mapsto\left[(\gamma_2(s)-\gamma_1(s))\left(x,\frac{|x|^2}{2},\frac{|x|^2}{2}\right)+(0,\gamma_1(s),\gamma_2(s))\right]\in\c\h^n,
\]
where $[(\gamma_1(s),\gamma_2(s))]$ are certain curves in $\c\h^1$ and 
$x\in\s^{n-1}\mapsto[x]\in\c\p^{n-1}$, $x\in\r\h^{n-1}\mapsto[x]\in\c\h^{n-1}$ and $x\in\r^{n-1}\mapsto x\in\c^{n-1}$ are the totally geodesic Lagrangian submanifolds in the  $(n\!-\!1)$-dimensional complex models.

The idea for constructing new minimal Lagrangian submanifolds in $\c\h^n$ is using the same curves as above and take, instead of these totally geodesic Lagrangian submanifolds, any minimal Lagrangian submanifold in the $(n-1)$-dimensional complex  models. In fact, it is straightforward to prove the following result.

\begin{proposition} \ 

\begin{description}

\item[a)] 

Given a solution $r(s)$ of the equation (1) (see Theorem 1) and a minimal Lagrangian immersion $\phi:N^{n-1}\rightarrow \c\p^{n-1}$ of a simply connected manifold $N$, $\Phi:\r\times N\rightarrow\c\h^n$ defined by
\[
\Phi (s,x) = \left[\left( 
\sinh r(s) \, e^{i \int_0^s \frac{dt}{\sinh^{n+1} r(t)} } \, \tilde{\phi}(x)  \, , \, 
\cosh r(s) \, e^{i \int_0^s \frac{\tanh^2 r(t) \, dt} {\sinh^{n+1} r(t)} }
\right)\right],
\]
is a minimal Lagrangian immersion in $\c\h^n$, where $\tilde{\phi}:N\rightarrow\s^{2n-1}$ is the horizontal lift of $\phi$ with respect to the Hopf fibration $\Pi:\s^{2n-1}\rightarrow\c\p^{n-1}$.
\item[b)] 

Given a solution $r(s)$ of the equation (3) (see Theorem 2) and a minimal Lagrangian immersion $\psi:N^{n-1}\rightarrow\c\h^{n-1}$ of a simply connected manifold N, $\Psi:\r\times N\rightarrow\c\h^n$ defined by
\[
\Psi (s,x) = \left[\left( 
\sinh r(s) \, e^{i \int_0^s \frac{\coth^2 r(t) \, dt}{\cosh^{n+1} r(t)} } \,  ,  \,
\cosh r(s) \, e^{i \int_0^s \frac{dt} {\cosh^{n+1} r(t)} } \, \tilde{\psi}(x)
\right)\right],
\]
is a minimal Lagrangian immersion in $\c\h^n$, where $\tilde{\psi}:N\rightarrow\h^{2n-1}_1$ is the horizontal lift of $\phi$ with respect to the Hopf fibration $\Pi:\h^{2n-1}\rightarrow\c\h^{n-1}$.

\item[c)] Given $\rho>0$, a minimal Lagrangian immersion $\eta:N^{n-1}\rightarrow\c^{n-1}$ of a simply-connected manifold $N$
and $f:N\longrightarrow \c$ satisfying $\Re f = |\eta|^2$ and $v(\Im f)=2\langle \eta_* v,J\eta \rangle $, for any vector $v$ tangent to $N$, the map $\Upsilon:\r\times N\rightarrow\c\h^n$ defined by
\[
\begin{array}{c}
\Upsilon(s,x) =
 \left[ e^{iA_{n+1}(s)} \left(
r(s)\eta(x),
\frac{1+r(s)^2(f(x)-1-2i A_{n+3}(s))}{2r(s)},
\right. \right. \\
\left. \left.
\frac{1+r(s)^2(f(x)+1-2i A_{n+3}(s))}{2r(s)}
\right) \right],
\end{array}
\]
where $A_n(s)=\int_0^s dt/r(t)^n $,
\(
r(s)=\rho \cosh^{\frac{1}{n+1}} ((n\!+\!1)s),
\)
is a minimal Lagrangian immersion in $\c\h^n$.
\end{description}
\end{proposition}

We observe that if we take $\eta:\r^{n-1}\rightarrow\c^{n-1}$ in Proposition 3 c) as the totally geodesic immersion $\eta (x)=x$, then $\Im f $ is constant and it is easy to prove that the corresponding immersion is congruent to the one given in Theorem 3.

It is interesting to note that the totally geodesic Lagrangian submanifolds of $\c\h^n$ can be also described in a similar way to the examples given in Theorems 1, 2 and 3. In fact, we can give three different descriptions of the totally geodesic Lagrangian submanifolds of $\c\h^n$:
\[
(s,x)\in\r^+\times\s^{n-1}\mapsto[(\sinh s\,x,\cosh s)]\in\c\h^n,
\]
\[
(s,x)\in\r\times\r\h^{n-1}\mapsto [(\sinh s,\cosh s\,x)]\in\c\h^n,
\]
and
\[
(s,x)\in\r\times\r^{n-1}\mapsto \left[e^s\left(x,\frac{|x|^2}{2},\frac{|x|^2}{2}\right)+(0,-\sinh s,\cosh s)\right]\in\c\h^n.
\]
In the three cases, the used curve is the geodesic $[(\sinh\,s,\cosh\,s)]$ of $\c\h^1$ passing through the point $[(0,1)]$. In the same way that in Proposition 3, we can construct new examples of minimal Lagrangian submanifolds of $\c\h^n$ following this idea.

\begin{proposition} \

\begin{description}

\item[a)] 

Given a minimal Lagrangian immersion $\phi:N^{n-1}\rightarrow \c\p^{n-1}$ of a simply connected manifold $N$, 
\[
\begin{array}{c}
\Phi:\r^+\times N\rightarrow\c\h^n \\
(s,x)\mapsto \left[\left( 
\sinh s \, \tilde{\phi}(x)  \, , \, 
\cosh s \right)\right],
\end{array}
\]
is a minimal Lagrangian immersion in $\c\h^n$, where $\tilde{\phi}:N\rightarrow\s^{2n-1}$ is the horizontal lift of $\phi$ with respect to the Hopf fibration $\Pi:\s^{2n-1}\rightarrow\c\p^{n-1}$. 
\item[b)] 

Given a minimal Lagrangian immersion $\psi:N^{n-1}\rightarrow\c\h^{n-1}$, 
\[
\begin{array}{c}
\Psi:\r\times N\rightarrow\c\h^n\\
 (s,x)\mapsto \left[\left( 
\sinh s \,  ,  \,
\cosh s \, \tilde{\psi}(x)\right)\right],
\end{array}
\]
is a minimal Lagrangian immersion in $\c\h^n$, where $\tilde{\psi}:N\rightarrow\h^{2n-1}_1$ is the horizontal lift of $\phi$ with respect to the Hopf fibration $\Pi:\h^{2n-1}\rightarrow\c\h^{n-1}$.
\item[c)] Given a minimal Lagrangian immersion $\eta:N^{n-1}\rightarrow\c^{n-1}$ of a simply-connected manifold $N$
and $f:N\longrightarrow \c$ satisfying $\Re f = |\eta|^2$ and $v(\Im f)=2\langle \eta_* v,J\eta \rangle $, for any vector $v$ tangent to $N$, 
\[
\begin{array}{c}
\Upsilon: \r \times N \rightarrow \c\h^n \\
(s,x)\mapsto \left[ 
e^s\left (\eta(x),f(x),f(x) \right )+(0,-\sinh s,\cosh s)\right ]
\end{array}
\]
is a minimal Lagrangian immersion in $\c\h^n$.
\end{description}
\end{proposition}

The examples described in Propositions 3 and 4 are unique in the following sense.
\begin{proposition}
Let $\gamma=(\gamma_1,\gamma_2):I\rightarrow \h^3_1$ a Legendre curve.
\begin{description}
\item[a)] Given a Lagrangian immersion $\phi:N^{n-1}\rightarrow \c\p^{n-1}$ of a simply connected manifold $N$, the map $\Phi:I\times N\rightarrow\c\h^n$ defined by
\[
\Phi(s,x)=\left [\left (\gamma_1(s)\tilde{\phi}(x),\gamma_2(s)\right )\right ],
\]
where $\tilde{\phi}:N\rightarrow\s^{2n-1}$ is a horizontal lift of $\phi$ with respect to the Hopf fibration, is a minimal Lagrangian immersion if and only if $\Phi$ is congruent to some of the examples given in Propositions 3,a) and 4,a).

\item[b)] Given a Lagrangian immersion $\psi:N^{n-1}\rightarrow \c\h^{n-1}$ of a simply connected manifold $N$, the map $\Psi:I\times N\rightarrow\c\h^n$ defined by
\[
\Psi(s,x)=\left [\left (\gamma_1(s),\gamma_2(s)\tilde{\psi}(x)\right )\right ],
\]
where $\tilde{\psi}:N\rightarrow\h^{2n-1}_1$ is a horizontal lift of $\psi$ with respect to the Hopf fibration, is a minimal Lagrangian immersion if and only if $\Psi$ is congruent to some of the examples given in Propositions 3,b) and 4,b).

\item[c)] Given a Lagrangian immersion $\eta:N^{n-1}\rightarrow \c^{n-1}$ of a simply-connected manifold $N$ and
$f:N\longrightarrow \c$ satisfying $\Re f = |\eta|^2$ and $v(\Im f)=2\langle \eta_* v,J\eta \rangle $, for any vector $v$ tangent to $N$,
the map $\Upsilon:I\times N\rightarrow\c\h^n$ defined by
\[
\Upsilon(s,x)=\left [(\gamma_2(s)-\gamma_1(s))\left (\eta(x),\frac{f(x)}{2},\frac{f(x)}{2}\right )+(0,\gamma_1(s),\gamma_2(s))  \right ]
\]
is a minimal Lagrangian immersion if and only if $\Upsilon$ is congruent to some of the examples given in Propositions 3,c) and 4,c).

\end{description}
\end{proposition}
{\it Proof : \/} 
In order to illustrate the idea of the proof, we only prove a).

Thanks to the properties of $\gamma$ and $\phi$, $\Phi$ is always a Lagrangian immersion. After a very long but straightforward computation,
we arrive at the horizontal lift $H^*$ of the
mean curvature $H$ of our Lagrangian immersion $\Phi$, which is given by 
$nH^*=a(s)J\tilde{\phi}_s + (n-1)(\gamma_1 H_\phi^*,0)/|\gamma_1|^2$,
where
\[ a=\frac{\langle \gamma'',J\gamma' \rangle}{|\gamma'|^4}+(n-1)
\frac{\langle \gamma_1',J\gamma_1 \rangle}{|\gamma_1|^2 |\gamma'|^2}. \]
If we suppose that $\Phi $ is minimal, necessarily $\phi $ is too since 
$H_\phi$ must be zero and, in addition, $a\equiv 0$.
We use this last equation  writing $\gamma $ as we did when proving Theorem 1. A similar reasoning
leads to the only two possible solutions for $r(s)$ corresponding to the solution of equation (1) or to the trivial solution $r(s)=s$.
$_{\diamondsuit}$

\begin{remark}
{\rm It is interesting to remark the parallelism between the constructions of minimal Lagrangian submanifolds of $\c\h^n$ above and the ones given in [CU2], [H] and [J1] when the ambient space is complex Euclidean space $\c^n$. In fact we can summarize some results given in the above papers to get the following result.

\vspace{0.2cm}

{\sc Proposition A} [CU2,H,J1]
{\em Let $\gamma:I\rightarrow\c^*$ be a regular curve and $\phi:N^{n-1}\rightarrow\c\p^{n-1}$ a Lagrangian immersion of a simply-connected manifold. Then $\Phi:I\times N\rightarrow\c^n$ defined by
\[
\Phi(s,x)=\gamma(s)\tilde{\phi}(x),
\]
where $\tilde{\phi}:N\rightarrow\s^{2n-1}$ is a horizontal lift of $\phi$ with respect to the Hopf fibration, is a minimal Lagrangian submanifold if and only if $\phi$ is minimal and $\gamma^n$ has curvature zero. }

\vspace{0.2cm}

 Then, up to rotations in $\c^n$, the curve $\gamma^n$ can be taken as $\gamma^n(s)=(s,c)$ with $c\geq 0$, (i.e.  $\Im\,\gamma^n=c$, where $\Im $ means imaginary part). So, up to dilations, there are only two possibilities: $c=0$ and $c=1$. In the first case, the examples constructed in this way are cones with links $\tilde{\phi}$, and in the second case the examples given in [CU2, Remark 1], [H, Theorem A] and [J1, Theorem 6.4].

The idea developed in Proposition 5 and Proposition A allows to construct a wide family of Lagrangian submanifolds, not necessarily minimal. This class of Lagrangian submanifolds has been deeply studied in [RU] when the ambient space is $\c^n$ and in [CMU] when the ambient space is $\c\p^n$ and $\c\h^n$. Among other things, it can be characterized by the existence of a closed and conformal vector field $X$ on the Lagrangian submanifold satisfying $\sigma(X,X)=\rho\,JX$, for a certain function $\rho$.
}
\end{remark}

\section{A geometric characterization}
As we pointed out in section 3.1, the examples of minimal Lagrangian submanifolds of $\c\h^n$ given in Theorem 1 are foliated by $(n-1)$-geodesic spheres of $\c\h^n$ centered at the point $[(0,\dots,0,1)]$. In a similar way it can be checked that $\{\Psi_{\rho}(\{s\}\times\r\h^{n-1})\,,s\in\r\}$ defines a foliation on the minimal Lagrangian submanifolds given in Theorem 2 by $(n-1)$-tubes over hyperplanes. Finally, $\{\Upsilon_{\rho}(\{s\}\times\r^{n-1})\,,s\in\r\}$ defines a foliation on the minimal Lagrangian submanifolds given in Theorem 3 by $(n-1)$-horospheres. In the following result we prove that the examples described in Theorems 1, 2 and 3 are the only admitting this kind of foliations.
\begin{theorem}
Let $\phi:M\rightarrow\c\h^n$ be a minimal Lagrangian immersion in $\c\h^n$. 
\begin{description}
\item[a)] If $\phi$ is foliated by (n-1)-geodesic spheres of $\c\h^n$, then either $\phi$ is totally geodesic or is congruent to an open subset of one of the examples described in Theorem 1.
\item[b)] If $\phi$ is foliated by (n-1)-tubes over hyperplanes of $\c\h^n$, then either $\phi$ is totally geodesic or is congruent to an open subset of one of the examples described in Theorem 2.
\item[c)] If $\phi$ is foliated by (n-1)-horospheres of $\c\h^n$, then either $\phi$ is totally geodesic or is congruent to an open subset of one of the examples described in Theorem 3.
\end{description}
\end{theorem}
{\it Proof of a):\/} Our submanifold $M$ is locally $I\times \s^{n-1}$, where $I$ is an interval of $\r$ with $0\in I$, and for each $s\in I$, $\phi\left (\{s\}\times \s^{n-1}\right )$ is an $(n-1)$-geodesic sphere of $\c\h^n$. So (see section 2.2) there exists a Lagrangian subspace
$$
\r\h^n_s = \{ [zX(s)]\in \c\h^n,z\in \c^{n+1},z=\bar{z} \}.
$$
where $X(s)\in U^1(n+1)$, and there exists $Y(s)\in SO^1_0(n+1)$, such that
\[
\phi(s,x)=\left [\left (\sinh r(s)\,x,\cosh r(s)\right )Y(s)X(s)\right ].
\]
Calling $A(s)=X(s)Y(s)$, we finally get
\[
\phi(s,x)=\left [\left (\sinh r(s)\,x,\cosh r(s)\right )A(s)\right ],
\]
with $A(s)\in U^1(n+1)$.
Then, $[(0,\dots,0,1)A(s)]$ and $r(s)$ are the center and the radius of the $(n-1)$-geodesic sphere $\phi(\{s\}\times \s^{n-1})$. 

If we denote
\[
\hat{\phi}(s,x)=\left (\sinh r(s)\,x,\cosh r(s)\right )A(s),
\]
then $\hat{\phi}$ is a lift (not necessarily horizontal) of $\phi$ to $\h^{2n+1}_1$. But (locally) Lagrangian immersions in $\c\h^n$ have horizontal lifts to $\h^{2n+1}_1$; so there exists a smooth function $\theta(s,x)$ such that $\tilde{\phi}=e^{i\theta}\hat{\phi}$ is a horizontal lift of $\phi$ to $\h^{2n+1}$. In particular,
\(
(d\tilde{\phi}_{(s,x)}(0,v),\tilde{\phi}(s,x))=0
\)
for any $v\in T_x\s^{n-1}$, which means that $d\theta(v)=0$ and so $\theta(s,x)=\theta(s)$. So our horizontal lift is given by
\[
\tilde{\phi}(s,x)=\left(\sinh r(s)\,x,\cosh r(s)\right)B(s),
\]
where $B(s)=e^{i\theta(s)}A(s)$. Moreover, as 
\(
\left( \tilde{\phi}_s,\tilde{\phi}\right )=0,
\)
then we obtain, using that $B(s)\in U^1(n+1)$, that
\[
(\sinh r(s)\,x,\cosh r(s))B'(s)S\bar{B}^t(s)(\sinh r(s)x,\cosh r(s))^t=0,
\]
for any $x\in\s^{n-1}$.

From $B(s)S\bar{B}^t(s)=S$ we deduce
$B'(s)S\bar{B}^t(s)+B(s)S\bar{B'}^t(s)=0$. So $B'(s)S\bar{B}^t(s)=V(s)+iU(s)$ where $V(s)$ and $U(s)$ are real matrixes with $V(s)+V(s)^t=0$ and $U(s)=U(s)^t$.
So the last equation becomes in
\[
(\sinh r(s)\,x,\cosh r(s))U(s)(\sinh r(s)x,\cosh r(s))^t=0,
\]
for any $x\in\s^{n-1}$. From this equation it is easy to obtain a smooth function $a(s)$ such that the matrix $U(s)$ is written as
\[
U(s)=a(s)\left( \begin{array}{c|c}
 I_n &  \\
\hline
& -\tanh^2r(s)
\end{array}\right),
\]
for any $s\in I$.

Now, we write $V(s) $ in the following way:
\[
V(s)=\left( \begin{array}{c|c}
 V_0(s) & -v^t(s) \\
\hline
v(s)& 0
\end{array}\right).
\]
Let $Z(s)$ be the solution to the following equation:
\[
Z'(s)+Z(s)V_0(s)=0,\quad Z(0)=I_n.
\]
As $V_0(s)+V_0^t(s)=0$, then $(Z(s)Z^t(s))'=0$ and so $Z(s)Z^t(s)=Z(0)Z^t(0)=I_n$. Then $Z(s)$ is a curve in $O(n)$ and we can reparametrize our immersion by
\[
(s,x)\in I\times\s^{n-1}\mapsto (s,xZ(s))\in I\times\s^{n-1},
\]
obtaining that
\[
\tilde{\phi}(s,x)=\left(\sinh r(s)\,x,\cosh r(s)\right)C(s),
\]
where 
\[
C(s)=\left( \begin{array}{c|c}
 Z(s) &  \\
\hline
& 1
\end{array}\right)B(s).
\]
Now, $C'(s)S\bar{C}^t(s)=W(s)+iU(s)$, where
\[
W(s)=\left( \begin{array}{c|c}
 0 & -w(s)^t \\
\hline
w(s) & 0
\end{array}\right),
\]
with $w(s)=v(s)Z^t(s)$.
\vspace{0.1cm}

Now we are going to use the minimality of our immersion. In order to do so, we first look for an orthonormal basis in our submanifold.
For any $x\in\s^{n-1}$, the vectors 
\[
z(s)=\tanh^{-1} r(s)\,w(s)-(\tanh^{-1} r(s)\,w(s)x^t)x \in \r^n
\]
are in $T_x\s^{n-1}$. Now, it is easy to check that $(1,-z(s))$, for any $\s\in I$, is a tangent vector to $M$ orthogonal to $(0,v)$, for any $v\in T_x\s^{n-1}$. 
Thus, an orthonormal basis of the submanifold $M=I\times\s^{n-1}$ at the point $(s,x)$ is
\[
e_1=\frac{(1,-z(s))}{|(1,-z(s))|}\,,\,e_i=\frac{(0,v_i)}{\sinh r(s)}\,,\,i=2,\dots,n,
\]
with $\{v_2,\dots,v_n\}$ an orthonormal basis of $T_x\s^{n-1}$. As $H=0$, in particular
\[
\langle \sum_{i=1}^n\sigma(e_i,e_i),J\tilde{\phi}_*(0,v)\rangle=0,
\]
for any $v\in T_x\s^{n-1}$. But it is easy to check that
\[
\langle \sigma(e_i,e_i),J\tilde{\phi}_*(0,v)\rangle=0,
\]
for $i=2,\dots,n$. So that the above equation becomes in
\[
\langle \sigma(e_1,e_1),J\tilde{\phi}_*(0,v)\rangle=0.
\]
Using the definition of $e_1$, we obtain that
\[
\langle \tilde{\phi}_{ss},J\tilde{\phi}_*(0,v)\rangle=2\langle (\tilde{\phi}_s)_*(0,z(s)),J\tilde{\phi}_*(0,v)\rangle,
\]
for any $v\in T_x\s^{n-1}$. Now from the properties of the second fundamental form of Lagrangian submanifolds, the definition of $z(s)$ and the fact $C'(s)=(W(s)+iU(s))SC(s)$, it is straightforward to prove that the last equation becomes in
\[
\frac{a(s)}{\cosh^2r(s)}w(s)v^t=0,
\]
for any $s\in I$, $v\in T_x\s^{n-1}$ and $x\in \s^{n-1}$. So $a(s)w(s)=0$ for any $s\in I$. 
So, if we define
\[
I_1=\{s\in I\,/\,a(s)=0\},\quad I_2=\{s\in I\,/\,w(s)=0\},
\]
we have that $I_1\cup I_2=I$.

First, we will work on the open set $I-I_2$, where $U(s)=0$ and so $C'(s)=W(s)SC(s)$. This implies that $C(s)$ are real matrixes and hence $C(s)\in O^1(n+1)$. As a consequence, $\phi((I-I_2)\times \s^{n-1})$ lies in $\r\h^n$ and so $\phi$ is totally geodesic on this open subset.

If we now work on the open set $I-I_1$, we have $W(s)=0$ and  then $C'(s)=iU(s)SC(s)$. Looking at the expression of $U(s)$, we can integrate the above equation obtaining that
\[
C(s)=\left( \begin{array}{c|c}

e^{i\int_{s_0}^s a(r)dr}I_n &  \\
\hline
& e^{-i\int_{s_0}^s a(r)\tanh^2rdr}
\end{array}\right).
\]
Therefore it is clear that, in this case, our immersion is invariant under the action of $SO(n)$ and then $\phi$, on the open set $(I-I_1)\times\s^{n-1}$, is one of the examples described in Theorem 1.

Finally, since the second fundamental forms of the examples given in Theorem 1 are non-trivial, using the connectedness of $I$, it cannot happen that Int$(I_1)\not=\emptyset$ and Int$(I_2)\not=\emptyset$. This finishes the proof of part a).
\vspace{0.2cm}

We omit the proof of b) because it is quite similar to the one given in a).
\vspace{0.2cm}

{\it Proof of c):\/} In this case our submanifold $M$ is locally $I\times\r^{n-1}$, where $I$ is an interval of $\r$ with $0\in I$, and for each $s\in I$, $\phi(\{s\}\times\r^{n-1})$ is an $(n-1)$-horosphere of some $\r\h^n_s$ embedded in $\c\h^n$ as a totally geodesic Lagrangian submanifold. So, following a similar reasoning like in the proof of a), we get 
$$
\phi(s,x)=\left [\hat{\phi}(s,x)\right ]=\left [f(x)A(s)\right ]
$$
where
\[
f(x)=\left (x,\frac{|x|^2}{2},\frac{|x|^2}{2}+1\right ).
\]
Then $\hat{\phi}$ is a lift (not necessarily horizontal) of $\phi$ to $\h^{2n+1}_1$. But (locally) Lagrangian immersions in $\c\h^n$ have horizontal lifts to $\h^{2n+1}_1$, so there exists a smooth function $\theta(s,x)$ such that $\tilde{\phi}=e^{i\theta}\hat{\phi}$ is a horizontal lift of $\phi$ to $\h^{2n+1}$. In particular,
\(
(d\tilde{\phi}_{(s,x)}(0,v),\tilde{\phi}(s,x))=0
\)
for any $v\in T_x\s^{n-1}$, which means that $d\theta(v)=0$ and so $\theta(s,x)=\theta(s)$. Thus, our horizontal lift is given by
\[
\tilde{\phi}(s,x)=f(x) B(s),
\]
where $B(s)=e^{i\theta(s)}A(s)$. Moreover, as 
\(
( \tilde{\phi}_s,\tilde{\phi} ) =0,
\)
we obtain that
\[
f(x)B'(s)S\bar{B}(s)^tf(x)^t=0,
\]
for any $s\in I$ and any $x\in\r^{n-1}$.
Again, a similar reasoning like in the proof of a) says that $B'(s)S\bar{B}(s)^t=V(s)+iU(s)$, where $V(s)$ and $U(s)$ are real matrixes with $V(s)+V(s)^t=0$ and $U(s)=U(s)^t$. So last equation becomes in
\[
f(x)U(s)f(x)^t=0,
\]
for any $s\in I$ and any $x\in\r^{n-1}$. From this equation it is easy to get that the matrix $U(s)$ is written as
\[
U(s)=a(s)\left(
\begin{array}{c|c|c}
I_{n-1} &   &  \\
\hline
 & 2 & -1 \\
\hline
 & -1 & 
\end{array} \right),
\]
for certain smooth function $a(s)$.

Now, we put
\[
V(s)=\left(
\begin{array}{c|c|c}
V_{0}(s) & -v_1(s)^t  & -v_2(s)^t \\
\hline
v_1(s) & 0 & -\rho(s) \\
\hline
v_2(s) & \rho(s) & 0
\end{array} \right)
\]
and let $Z(s)$ be the solution to the following differential equation
\[
Z'(s)+Z(s)V_0(s)=0,\quad Z(0)=I_n.
\]
From $V_0(s)+V_0^t(s)=0$, it follows that $(ZZ^t)'(s)=0$, and so $Z(s)Z^t(s)=Z(0)Z^t(0)=I_n$. This means that $Z(s)$ is a curve in $O(n\!-\!1)$. We can now reparametrize our immersion by
\[
(s,x)\in I\times\r^{n-1}\mapsto (s,xZ(s))\in I\times\r^{n-1},
\]
obtaining that
\[
\tilde{\phi}(s,x)=f(x)C(s),
\]
where 
\[
C(s)=\left( \begin{array}{c|c}
 Z(s) &  \\
\hline
& I_2
\end{array}\right)B(s).
\]
Now, it is easy to check that $C'(s)S\bar{C}^t(s)=W(s)+iU(s)$, where
\[
W(s)=\left(
\begin{array}{c|c|c}
 & -w_1^t(s)  & -w_2^t(s) \\
\hline
w_1(s) &  & -\rho (s) \\
\hline
 w_2(s)& \rho (s)& 
\end{array} \right)
\]
with $w_i(s)=v_i(s)Z^t(s),\, i=1,2$.
\vspace{0.1cm}

Now we are going to use the minimality of our immersion. In order to do so, first  we are going to find an orthonormal basis in our submanifold at $(s,0)$.
It is easy to check that $(1,-w_2(s))$ is a tangent vector to $M$ in $(s,0)$, orthogonal to $(0,v)$ for any $v\in\r^{n-1}$. 
So an orthonormal basis of the submanifold $M=I\times\r^{n-1}$ at the point $(s,0)$ is
\[
e_1=\frac{(1,-w_2(s))}{|(1,-w_2(s))|}\,,\,e_i=(0,v_i)\,,\,i=2,\dots,n,
\]
with $\{v_2,\dots,v_n\}$ an orthonormal basis of $\r^{n-1}$. As $\phi$ is a minimal immersion, in particular we have that 
\[
\langle H(s,0),J\tilde{\phi}_*(0,v)\rangle=0,
\]
for any $v\in\r^{n-1}$. But it is easy to check that
\[
\langle \sigma(e_i,e_i),J\tilde{\phi}_*(0,v)\rangle=0,
\]
for $i=2,\dots,n$. In this way the above equation becomes 
\[
\langle \sigma(e_1,e_1),J\tilde{\phi}_*(0,v)\rangle=0.
\]
Using the definition of $e_1$, we obtain that
\[
\langle \tilde{\phi}_{ss},J\tilde{\phi}_*(0,v)\rangle=2\langle (\tilde{\phi}_s)_*(0,w_2(s)),J\tilde{\phi}_*(0,v)\rangle,
\]
for any $v\in\r^{n-1}$. From the properties of the second fundamental form of Lagrangian submanifolds and the fact $C'(s)=(W(s)+iU(s))SC(s)$, it is straightforward to prove that the last equation becomes in
\[
a(s)(w_1(s)+w_2(s))v^t=0,
\]
for any $s\in I$, $v\in\r^{n-1}$. So $a(s)(w_1(s)+w_2(s))=0$ for any $s\in I$. 
So, if we define
\[
I_1=\{s\in I\,/\,a(s)=0\},\quad I_2=\{s\in I\,/\,w_1(s)+w_2(s)=0\},
\]
we have that $I_1\cup I_2=I$.

First, we will work on the open set $I-I_2$. There $U(s)=0$ and so $C'(s)=W(s)SC(s)$. This implies that $C(s)$ are real matrixes and hence $C(s)\in O^1(n+1)$. As a consequence, $\phi((I-I_2)\times \r^{n-1})$ lies in $\r\h^n$ and so $\phi$ is totally geodesic on this open subset.

If we now work on the open set $I-I_1$, we have that 
\[
W(s)=\left(
\begin{array}{c|c|c}
 & -w_1^t(s)  & w_1^t(s) \\
\hline
w_1(s) &  & -\rho (s) \\
\hline
 -w_1(s)& \rho (s) & 
\end{array} \right).
\]
If $w(s)$ is a solution of $w'(s)+\rho(s)w(s)-w_1(s)=0$, we can reparametrize our immersion as
\[
(s,x)\in (I-I_1)\times\r^{n-1}\mapsto (s,x+w(s))\in(I-I_1)\times\r^{n-1},
\]
so that the immersion is given by
\[
\tilde{\phi}(s,x)=f(x)D(s),\]
where 
\[
D(s)=\left(
\begin{array}{c|c|c}
 I_{n-1}& w^t(s)  & w^t(s) \\
\hline
-w(s) &1-\lambda  & -\lambda \\
\hline
 w(s)& \lambda & 1+\lambda
\end{array} \right)C(s),
\]
where $\lambda=|w(s)|^2/2$.
Now it is easy to check that $D'(s)S\bar{D}^t(s)=W^1(s)+iU(s)$ with
\[
W^1(s)=\left(
\begin{array}{c|c|c}
 &   &  \\
\hline
 &  & -\rho (s) \\
\hline
 & \rho (s) & 
\end{array} \right).
\]
We consider the matrixes
\[
Y(s)=\left(
\begin{array}{c|c|c}
I_{n-1} &   &  \\
\hline
 & \cosh \int \rho(s) &  -\sinh \int \rho(s)\\
\hline
 & -\sinh \int \rho(s) & \cosh \int \rho(s)
\end{array} \right)
\]
and define $F(s)=Y(s)D(s)$. Note that 
\[
Y(s)^{-1}=\left(
\begin{array}{c|c|c}
I_{n-1} &   &  \\
\hline
 & \cosh \int \rho(s) &  \sinh \int \rho(s)\\
\hline
 & \sinh \int \rho(s) & \cosh \int \rho(s)
\end{array} \right).
\]
Then it is easy to get that $F'(s)S\bar{F}^t(s)= i F(s)U(s)\bar{F}^t(s)$ and from here we arrive at the linear differential equation
$F'(s)=G(s)F(s)$, where
\[
G(s)= i a(s) \left(
\begin{array}{c|c|c}
I_{n-1} &   &  \\
\hline
 & 1+ \cosh 2\int \rho(s) &  1+\sinh 2\int \rho(s)\\
\hline
 & -(1+\sinh 2\int \rho(s)) & 1-\cosh 2\int \rho(s)
\end{array} \right),
\]
whose solution can be written as $F(s)=e^{\int G(s)}$. Therefore $D(s)=e^{\int G(s)}Y(s)^{-1}$
and now it can be easily checked that the immersion is invariant under the action of $SO(n\hspace{-3pt}-\hspace{-3pt}1)\propto\r^{n-1}$ and so $\phi$, on $(I-I_1)\times\r^{n-1}$, is one of the examples given in Theorem 3. Finally, as the second fundamental forms of the examples given in Theorem 3 are non trivial, using the connectedness of $I$ it cannot happen that  Int$(I_1)\not=\emptyset$ and Int$(I_2)\not=\emptyset$. This finishes the proof of part c).$_{\diamondsuit}$

\section{Minimal Lagrangian submanifolds in $\c\p^n$}

As we mentioned in the introduction, in this section we are going to describe (without proofs) the corresponding results when the ambient space is the complex projective space $\c\p^n$. 

If $U(n+1)$ is the unitary group of order $n+1$, then $PU(n+1)=U(n+1)/\s^1$ is the group of holomorphic isometries of $(\c\p^n,\langle,\rangle)$. We consider the special orthogonal group $SO(n)$ acting on $\c\p^n$ as a subgroup of holomorphic isometries in the following way:
\[
A\in SO(n)\longmapsto \left[\left(\begin{array}{c|c}

A &  \\
\hline
&  1
\end{array}\right)\right]\in PU(n+1),
\]
where $[\;\;]$ stands for class in $U(n+1)/\s^1$.

The unit sphere $\s^n$ can be isometrically immersed in $\c\p^n$ as a totally geodesic Lagrangian submanifold in the standard way
\[ 
x\in\s^n\mapsto [x]\in\c\p^n.
\]
This immersion projects in the totally geodesic Lagrangian embedding of the real projective space $\r\p^n$ in $\c\p^n$. Moreover, up to congruences, it is the only totally geodesic Lagrangian submanifold of $\c\p^n$. It is interesting to note that the totally umbilical submanifolds of $\c\p^n$ (which were classified in [ChO]) are either totally geodesic or umbilical submanifolds of totally geodesic Lagrangian submanifolds. So, {\em up to conguences, the $(n-1)$-dimensional totally umbilical (non-totally geodesic) submanifolds of $\c\p^n$ are the umbilical hypersurfaces of $\r\p^n$ embedded in $\c\p^n$ in the above way}. In this case, the umbilical hypersurfaces of $\r\p^n$ are the geodesic spheres. We will refer to these examples as {\em $(n-1)$-geodesic spheres of $\c\p^n$}.

\begin{theorem}
Let $\phi:M\rightarrow\c\p^n$ be a minimal (non-totally geodesic) Lagrangian immersion.
\begin{description}
\item[a)] $\phi$ is invariant under the action of $SO(n)$ if and only if $\phi$ is locally congruent to one of the immersions in the following 1-parameter family of minimal Lagrangian immersions $\{\Phi_{\rho}:\r\times\s^{n-1}\rightarrow\c\p^n\,/\rho\in\r^+\}$, given by 
\[ 
\Phi_{\rho} (s,x) = \left[\left( 
\sin r(s) \, e^{-i \int_0^s \frac{dt}{\sin^{n+1} r(t)} } \, x  \, , \, 
\cos r(s) \, e^{i \int_0^s \frac{\tan^2 r(t) \, dt} {\sin^{n+1} r(t)} }
\right)\right],
\]
where $r(s)$, $s\in \r$, is the only solution to
\begin{equation}
r''\sin r \cos r  = (1-(r')^2)(n\cos^2 r-\sin^2 r), \, r(0)=\rho, \,
r'(0)=0.
\end{equation}
\item[b)] $\phi$ is foliated by (n-1)-geodesic spheres of $\c\p^n$ if and only if $\phi$ is locally congruent to one of the examples described in a).
\end{description}
\end{theorem}

\begin{remark}
{\rm In this case, $r(s)=\arctan\sqrt n$ gives a constant solution to equation (5). The corresponding minimal Lagrangian immersion $\Phi:\r\times\s^{n-1}\rightarrow\c\p^n$ is given by
\[
\Phi(s,x)=\left[ \frac{1}{\sqrt{n+1}}\,\left(\sqrt {n}e^{-is}\,x,e^{ins}\right)\right],
\]
which provides a minimal Lagrangian immersion $\Phi:\s^1\times\s^{n-1}\rightarrow\c\p^n$ defined by
\[
\Phi(e^{it},x)=\left[ \frac{1}{\sqrt{n+1}}\left(\sqrt {n}\,e^{\frac{-it}{n+1}}\,x,e^{\frac{int}{n+1}}\right)\right].
\]
If $h:\s^1\times\s^{n-1}\rightarrow\s^1\times\s^{n-1}$ is the diffeomorphism $h(e^{it},x)=(-e^{it},-x)$, then $\Phi$ induces a minimal Lagrangian embedding $(\s^1\times\s^{n-1})/h\rightarrow\c\p^n$, which is a very well-known example studied by Naitoh ([N], Lemma 6.2).}
\end{remark}

\begin{remark}
{\rm By studying the energy integral of equation (5) given by
\[
(r')^2+\frac{\sin^{2n}\rho\cos^2\rho}{\sin^{2n}r\cos^2r}=1,
\]
it is easy to check (when $r$ is not the constant solution) that the orbits $s\mapsto(r(s),r'(s))$ are closed curves. Hence, all the solutions of equation (5) are periodic functions. However, not all the corresponding minimal Lagrangian submanifolds are embedded. In fact, in [CU1] minimal Lagrangian surfaces invariant by a $1$-parameter group of holomorphic isometries of $\c\p^2$ were classified, obtaining the examples given in Theorem 5 with $n=2$ as a particular case. As the solutions of (5) for $n=2$ are elliptic functions (see [CU1]), it is not difficult to check that, except the Clifford torus, the examples given there do not provide embedded minimal Lagrangian tori. It may be interesting to point out here that recently Goldstein [G] has constructed minimal Lagrangian tori in Einstein-Kaehler manifolds with positive scalar curvature.}
\end{remark}

Now we give a method to produce examples of minimal Lagrangian submanifolds of $\c\p^n$.
\begin{proposition}
Let $\phi:N^{n-1}\rightarrow \c\p^{n-1}$ be a minimal Lagrangian immersion of a simply connected manifold $N$, and $\tilde{\phi}:N\rightarrow\s^{2n-1}$ the horizontal lift of $\phi$ with respect to the Hopf fibration $\Pi:\s^{2n-1}\rightarrow\c\p^{n-1}$.
\begin{description}
\item[a)]
Given a solution $r(s)$ of the equation (5) (see Theorem 5),  $\Phi:\r\times N\rightarrow\c\p^n$ defined by
\[
\Phi (s,x) = \left[\left( 
\sin r(s) \, e^{-i \int_0^s \frac{dt}{\sin^{n+1} r(t)} } \, \tilde{\phi}(x)  \, , \, 
\cos r(s) \, e^{i \int_0^s \frac{\tan^2 r(t) \, dt} {\sin^{n+1} r(t)} }
\right)\right]
\]
is a minimal Lagrangian immersion in $\c\p^n$.
\item[b)] The map 
\[
\begin{array}{c}
\Phi:(0,\frac{\pi}{2})\times N\rightarrow\c\p^n\\
 \\
(s,x)\mapsto \left[\left( 
\sin s \, \tilde{\phi}(x)  \, , \, 
\cos s \right)\right]
\end{array}
\]
is a minimal Lagrangian immersion in $\c\p^n$. 
\item[c)] Let $\gamma=(\gamma_1,\gamma_2):I\rightarrow \s^3$ be a Legendre curve. The map $\Phi:I\times N\rightarrow\c\p^n$ defined by
\[
\Phi(s,x)=\left [\left (\gamma_1(s)\tilde{\phi}(x),\gamma_2(s)\right )\right ]
\]
is a minimal Lagrangian immersion if and only if $\Phi$ is congruent to some of the examples given in a) and b).
\end{description}
\end{proposition}

\vspace{1cm}

\begin{tabular}{ll}
{\sc addresses}: &  \\
(first and second authors) & (third author) \\
Departamento de Matem\'{a}ticas & Departamento de Geometr\'{\i}a  \\
Escuela Polit\'{e}cnica Superior & y Topolog\'{\i}a \\
Universidad de Ja\'{e}n & Universidad de Granada\\
23071 Ja\'{e}n & 18071 Granada \\
SPAIN & SPAIN \\
 & \\
{\sc e-mails}: & \\
(first author) {\tt icastro@ujaen.es} \\
(second author)  {\tt crodri@ujaen.es} \\
(third author) {\tt furbano@goliat.ugr.es} \\
\end{tabular}


\begin{thebibliography}{[BEM1]}

\bibitem [{\bf CMU}]{}
I. Castro, C.R. Montealegre, F. Urbano, {\em Closed conformal vector fields and Lagrangian submanifolds in complex space forms}, Pacific J. Math. (to appear).

\bibitem [{\bf CU1}]{}
I. Castro, F. Urbano, {\em New examples of minimal Lagrangian tori in the complex projective plane}, Manuscripta Math. {\bf 85} (1994), 265--281.

\bibitem [{\bf CU2}]{}
I. Castro, F. Urbano, {\em On a minimal Lagrangian submanifold of $C^n$ foliated by spheres},
Michigan Math. J. {\bf 46} (1999), 71--82.

\bibitem [{\bf ChO}]{}
B.Y. Chen, K. Ogiue, {\em Two theorems on Kaehler manifolds}, Michigan Math. J. {\bf 21} (1974), 225--229.

\bibitem [{\bf G}]{}
E. Goldstein, {\em Minimal Lagrangian tori in Kaehler--Einstein manifolds}, math.DG/0007135, 2000.

\bibitem [{\bf HL}]{}
R. Harvey, H.B. Lawson, {\em Calibrated geometries}, Acta Math. {\bf 148} (1982), 47--157.

\bibitem [{\bf H}]{}
M. Haskins, {\em Special Lagrangian cones}, math.DG/0005164, 2000.

\bibitem [{\bf J1}]{}
D.D. Joyce, {\em Special Lagrangian m-folds in $\c^m$ with symmetries}, math.DG/0008021, 2000.

\bibitem [{\bf J2}]{}
D.D. Joyce, {\em Constructing special Lagrangian m-folds in $\c^m$ by involving quadrics}, math.DG/0008155, 2000. 

\bibitem [{\bf J3}]{}
D.D. Joyce, {\em Evolution equations for special Lagrangian 3-folds in $\c^3$}, math.DG/0010036, 2000.

\bibitem [{\bf J4}]{}
D.D. Joyce, {\em Ruled special Lagrangian 3-folds in $\c^3$}, math.DG/0012060, 2000.

\bibitem [{\bf N}]{}
H. Naitoh, {\em Isotropic submanifolds with parallel second fundamental form in $\p^m(\c)$}, Osaka J.Math. {\bf 18} (1981), 427--464.

\bibitem [{\bf O}]{}
Y.G. Oh, {\em Second variation and stabilities of minimal Lagrangian submanifolds}, Invent. Math. {\bf 101} (1990), 501--519.

\bibitem [{\bf RU}]{}
A. Ros, F. Urbano, {\em Lagrangian submanifolds of $\c^n$ with conformal Maslov form and the Whitney sphere}, J. Math. Soc. Japan {\bf 50} (1998), 203--226.

\bibitem [{\bf SYZ}]{}
A. Strominger, S.T. Yau and E. Zaslow, {\em Mirror symmetry is T-duality}, Nuclear Physics B {\bf 479} (1996), 243--259.

\end{thebibliography}
\end{document}